\documentclass[12pt]{amsart}
\usepackage{amsmath,amsthm,epsfig,amssymb}
\usepackage{amsthm,amsfonts,amssymb,amscd}
\usepackage[all]{xy}
\newtheorem{thm}[subsection]{Theorem}

\newtheorem{lem}[subsection]{Lemma}
\newtheorem{corol}[subsection]{Corollary}
\newtheorem{rem}[subsection]{Remark}
\theoremstyle{definition}
\newtheorem{Def}[subsection]{Definition}

\newtheorem{proposition-definition}[subsection]{Proposition-Definition}

\newcommand{\CC}{{\mathbb C}}

\newcommand{\PP}{{\mathbb P}}

\newcommand{\OOO}{{\mathcal O}}

\author{A. El Mazouni}

\address{Laboratoire de Math\'ematiques de Lens EA 2462
Facult\'e des Sciences Jean Perrin
Rue Jean Souvraz, SP18
F-62307 LENS  Cedex France}
\email{mazouni@euler.univ-artois.fr}

\author{F. Laytimi}

\address{F. L.: Math\'ematiques - b\^{a}t. M2, Universit\'e Lille 1,
F-59655 Villeneuve d'Ascq Cedex, France}
\email{fatima.laytimi@math.univ-lille1.fr}

\author{D.S. Nagaraj}

\address{Institute of Mathematical Sciences C.I.T. campus, Taramani, 
chennai 600113,India}
\email{dsn@imsc.res.in}

\title{Projections of Veronese surface and Morphisms from projective plane to Grassmannian}
\begin{document}

\begin{abstract}  

In this note we describe the image of $\PP^2$ in $ Gr(2, \CC^{4})$
under a morphism given by 
a rank two vector bundle on $\PP^2$ with Chern classes $(2,2).$ 

\end{abstract}

\maketitle

\section{Introduction} \setcounter{page}{1}
We denote by $\PP^2$ the projective plane over the field $\CC$ 
of complex numbers and 
by $ Gr(2, \CC^{4})$  the Grassmannian variety of two dimensional
quotient spaces of $ \CC^{4}.$

 Let $Q$ be a rank two vector bundle on $\PP^2$ generated by global sections. 
Then $Q$ can be generated by at most four linearly
independent section. Assume that $Q$ is generated by four linearly
independent sections but is not generated by less number of sections.
If we fix a set $S$  of four linearly independent 
global sections generating $Q$  we get a non constant morphism
$$ \phi_{S} : \PP^2  \to Gr(2, \CC^{4}). $$
The aim of this article  is to study the properties of the
image of $\PP^2$ under such a morphism $\phi_{S}.$ 

According to a result of Tango \cite{Ta}, if $\phi_{S}$ is an 
embedding then the pair of Chern classes $(c_1(Q), c_2(Q))$ is 
either $(2,1)$ or $ ((2,3).$ 

It is interesting to know what is the image of $\PP^2$ under
a general morphism $\phi_S : \PP^2 \to Gr(2 , \CC^{4}).$  For example, What are 
the possible singularities of the image ?  What are the 
defining equations of the image? 

In this article  we focus on the image of $\PP^2$ in 
$Gr(2, \CC^{4})$  associated to a  
globally generated vector bundle $Q$ with Chern class
 pair  $(c_1(Q), c_2(Q)) = (2,2).$
More precisely we have the following Theorem(\ref{image}):
\begin{thm}
 Let $\phi:\PP^2 \to Gr(2,\CC^{4})$ be a morphism. Assume that 
 $c_1(Q)=2$ and $c_2(Q)= 2,$ where $Q$ is the pull back
 by $\phi$ of the universal rank two quotient bundle on $Gr(2, \CC^{4}).$
 Then the image of $\PP^2$ in $Gr(2, \CC^{4})$ is 
 \begin{itemize}
  \item[a)] either a complete intersection of two independent hyperplanes,
  
  \item[b)] or is a complete intersection of a hyperplane and a quadric.
 \end{itemize}
Here a hyperplane(respectively, quadric)
 means a divisor in the class of the ample generator (respectively, twice the 
class of the  ample generator)
 of the the Picard group of $Gr(2, \CC^{4}).$ In case of a) the image is 
 isomorphic to a cone in $\PP^3$ over a conic. In the  case of b) the image 
  is isomorphic to base locus of a pencil consisting of  singular
quadrics in $\PP^4$ of rank $3$ and $4.$ More over the image surface is 
singular exactly along a line of $\PP^4$ contained in the surface.

\end{thm}

\section{Special Projections of Veronese Surface}
 
 Our study of the image of $\PP^2$ in $ Gr(2, \CC^{4})$ under the 
 morphism given by a rank two vector bundle on $\PP^2$ with Chern classes $(2,2),$ 
depends on the study of the projections of Veronese surface $V$ in $\PP^5$ 
from a special point not on the surface or a special line not intersecting the Veronese surface. Here special point means a point on the secant variety 
of $V$ and a special line means a line contained in the secant variety. 
In this section we recall some facts  about the Veronese surface 
(see, \cite{BCGM}, \cite{Hari} and \cite{FuHa} for details) and
deduce some results about special projections.

Veronese surface is the only non-degenerate (i.e., not contained in a
hyperplane)
non-singular surface  in $\PP^5$ which can be projected isomorphically
to $\PP^4.$  A general projection has this
property. It is well known that the Veronese surface is the unique closed orbit
for the natural action of the algebraic group $PGL(3, \CC^3)$ on $\PP^5.$
In fact if we identify $\PP^5$ with the space of conics in $\PP^2,$ there
are three orbits namely, 
\begin{itemize} 
 \item the set of all non-singular quadrics
 
 \item the set of all pairs of distinct lines
 
 \item the set of all double lines.
 
\end{itemize}
The set of double lines is the Veronese surface. 
(See, \cite[p.120-21]{Hari}). 

 The following remarks shed light on the image of Veronese surface under 
 special projections.
 \begin{rem}\label{point} Let $V$ be  a Veronese surface in
$P(H^0(\PP^2,\mathcal{O}_{\PP^2}(2)))$ and $sec(V)$ be 
its secant variety. Let $p \in sec(V)\setminus V$ be a point. Then the
the projection of $\PP^5\setminus \{p\}$ from $p$ to a hyperplane 
$H\simeq \PP^4$ 
in $\PP^5$ not containing $p$ maps $V$ onto a singular surface $V_p$
of degree four.  $V_p$ is cut out by intersection of quadrics 
of a linear pencil of quadrics in $\PP^4,$ and each qudaric in the 
pencil is either has rank three or four. More over $V_p$ is singular
along a line $L_p$ of $\PP^4.$ The projection map $f: V \to V_p$ induces an 
isomorphism from $V\setminus f^{-1}(L_p)$ to $V_p\setminus L_p,$ and
$f^{-1}(L_p) $ is a conic and the map $f$ restricted to $f^{-1}(L_p) $
to $L_p$ is a ramified covering of degree two.  
\end{rem}
See \cite[p. 366, 10.5.5.]{BCGM} for details. The only thing missing 
there is the statement that $V_p$ is cut out by intersection of quadrics
of $\PP^4,$ in a linear pencil of quadrics in $\PP^4,$ and each quadric in the 
pencil is either has rank three or four. This can be seen in several
ways. One of the way  is to use the fact that the set of  all pairs
of  distinct 
lines is one orbit for the action for $PGL(3, \CC^3)$ namely 
$sec(V)\setminus V$ and to prove that required property holds for the
projection {\em corresponds to  a particular point} 
$p \in sec(V)\setminus V.$  
The image morphism  $\PP^2 \to \PP^4$  given by  
 $$(x,y,z)\mapsto (z^2, -xy, -y^2+yz, x^2, xy-xz)$$ is one such projection.
 The image is cut out by the pencil of quadrics 
 $\lambda (Z_1Z_4-Z_2Z_3)+\mu ((Z_1+Z_4)^2-Z_0Z_3), $ 
 where $Z_i \, \, (0\leq i \leq 4)$ 
 are the homogeneous coordinate functions on $\PP^4.$ 
 
 \begin{rem}\label{line} Let $V$ be  a Veronese surface in
$P(H^0(\PP^2,\mathcal{O}_{\PP^2}(2)\simeq \PP^5$ and $sec(V)$ be 
its secant variety. Let $\ell \subset sec(V) $ be a line of $\PP^5$ 
such that
$V \cap \ell = \emptyset.$  Then the
the projection of $\PP^5\setminus \ell$ from $\ell$ to a linear subspace 
$ L\simeq \PP^3$ 
in $\PP^5$ not meeting $\ell$ maps $V$ onto a qudaric surface $V_{\ell},$
which is a cone over a quadric.   The projection map $f: V \to V_{\ell}$ 
is generically two to one. 
\end{rem}
 Existence of lines in $sec(V)$ with the above property follows from
 \cite[p. 361, 10.4.]{BCGM}. For example, for a fixed line 
 $L_0 \subset \PP^2$ the plane $[L_0](\subset sec(V))\subset \PP^5 $ 
 defined by 
 $$ [L_0]:= \{ L_0\cup M| M \subset \PP^2, \text{a line} \} $$
  meet  the Veronese surface $V$ in a single point 
  $l_0.$ Here we have identified $V$ as set of all lines in $\PP^2$
  under the map $L \mapsto L^2.$ Under this identification 
  $sec(V)$ corresponds to set of reducible conics. Hence any 
  line in the plane $[L_0]$ 
 not passing through $l_0$has the required property.

\section{Morphisms from $\PP^2$ to $Gr(2,\CC^{4}).$}

For a vector bundle $E$ on $\PP^2$
the bundle $E\otimes \OOO_{\PP^2}(k)$ is denoted by $E(k).$

A globally generated rank two vector bundle $Q$ on $\PP^2$
can be generated by at most four linearly independent global sections.
If we take a set of generators consisting of at most four sections of 
$Q$ we get  a morphism from 
$\PP^2$ to $ Gr(2, \CC^{4}).$ Note that the bundle $Q$ is generated by two linearly independent
sections if and only if $Q \simeq \OOO^2_{\PP^2}.$ This happens if and only
the morphism from $\PP^2$ to $Gr(2, \CC^{4})$ is  constant.
Moreover the bundle  $Q$ generated by three linearly independent sections
if and only if the morphism from $\PP^2$ to $ Gr(2, \CC^{4})$ factors through a linear
$\PP^2\simeq Gr(2, \CC^{3})$ contained in $ Gr(2, \CC^{4}).$

{\Def Let $Q$ be a rank two vector bundle on $\PP^2 $ generated by global sections. 
Assume that $Q$ cannot be generated by less than four independent
sections. If $S$ is a set of four independent global sections generating
$Q,$ then we get morphism
$$\phi_S :\PP^2 \to Gr(2, \CC^{4}).$$
We call such a morphism $\phi_S$ a non special morphism.
(Generally, we use the the notation $\phi$ instead of $\phi_S.$)

{\rem 1) Let $\phi_S : \PP^2 \to Gr(2, \CC^{4})$ be a non special morphism obtained from a
rank two vector bundle $Q.$ Then 
the pull back to $\PP^2$ by $\phi_S$ of the universal quotient bundle on 
$Gr(2, \CC^{4})$ is equal to $Q.$ Since the morphism
$\phi_S$ is non special  ${\rm det}(Q) = \OOO_{\PP^2}(d)$ for some 
$d > 0,$ i.e., $c_1(Q) > 0.$ As $Q$ is generated by sections we see that
$c_2(Q)\leq 0.$ A rank two bundle $Q$ generated by sections has $c_2(Q)=0$ 
implies $Q = \OOO_{\PP^2}\oplus \OOO_{\PP^2}(d).$ This implies 
$\phi_S : \PP^2 \to Gr(2, \CC^{4})$ is a special morphism, 
a contradiction to the assumption. Thus we must have $c_2(Q)> 0.$

2) Let $\phi_S : \PP^2 \to Gr(2, \CC^{4})$ be a non special morphism 
as above.  If $ p: Gr(2, \CC^{4}) \to \PP^5 $ is the Plucker imbedding, then
note that the morphism 
$$ p\circ\phi_{S} : \PP^2 \to \PP^5 $$ 
may not be non degenerate in the usual sense. In other words 
the image of $\PP^2$ in $\PP^5 $ under $ p\circ\phi_{S}$ may
be very well contained in a hyperplane of $\PP^5. $}

For the study of non special morphisms one need to know what are the
globally generated rank two vector bundles on $\PP^2.$  
 In our previous paper
\cite{ALN} we obtained some partial
results about the  possible Chern classes $(c_1(Q), c_2(Q))$ of rank two vector bundles $Q$
on $\PP^2$ generated by four sections. In \cite{PE} Ph.Ellia determined
the Chern classes of rank two globally generated vector bundles on $\PP^2.$
His  result gives the  complete numerical characterization of such bundles.
Note that a  rank two globally generated vector bundles on $\PP^2$ 
can be generated by $4$ sections and hence gives rise to a morphism
from $\PP^2$ to $ Gr(2, \CC^{4}). $
Globally generated vector bundle on projective spaces 
with special Chern classes are studied in \cite{CU1} and \cite{CU2}.

According to a result of Tango \cite{Ta}, if $\phi_S$ is a non special imbedding
then the Chern class pair  $(c_1(Q), c_2(Q)),$ of  $Q$ is either $(2,1)$ or $(2,3).$ 

The aim of this note is to investigate the properties of the image
of $\PP^2$ under  $\phi_S$ for the case $c_1(Q) =2,\, c_2(Q) =2.$  For such a
bundle we have the following:

\begin{lem}\label{lem0}
Let $Q$ be a rank two vector bundle 
on $\PP^2$ with $c_1(Q) =2,\, c_2(Q) =2.$ 
If $Q$ is generated by sections then $Q$ is semi-stable. 
\end{lem}

{\bf Proof:}  Assume that $Q$ is not semi-stable,
 then by Lemma (3.1) of \cite{Ha2}  $h^0(Q(-2)) \neq 0.$
 Let $k$ be the largest integer such that $h^0(Q(-k)) \neq 0.$
 Note that $k \geq 2$ and
 there is an exact sequence of sheaves
\begin{equation}\label{eq0}
 0 \to \OOO_{\PP^2} \to Q(-k) \to I_Z(-{2k}+2) \to 0,
 \end{equation}
where $I_Z$ is the ideal sheaf of a zero-dimensional closed sub-scheme 
$Z$ of length  $k^2-2k+2= (k-1)^2+1.$ Now tensoring the exact sequence
\eqref{eq0} with
the line bundle $\OOO_{\PP^2}(k)$ we get the following exact sequence
\begin{equation}\label{eq01}
 0 \to \OOO_{\PP^2}(k) \to Q \to I_Z(-k+2) \to.0
 \end{equation}
Since $Q$ generated by sections $k=2$ and $Z =\emptyset.$  This implies 
$c_2(Q) = 0, $ a contradiction. This contradiction proves the lemma. $\hfill{\Box}$

\begin{lem}\label{lem1}
Let $Q$ be a semi-stable rank two vector bundle 
on $\PP^2$ with $c_1(Q) =2,\, c_2(Q) =2.$ 
Then $Q$ is generated by four independent  sections and hence  
there is a surjective 
morphism  of bundles
$$ \OOO_{\PP^2}^4 \to Q \to 0$$
which determines a morphism $\phi : \PP^2 \to Gr(2, \CC^{4}).$
\end{lem}
 
 {\bf Proof:}  If $Q$ is a rank two vector bundle on $\PP^2$ with 
 $c_1(Q) =2,\, c_2(Q) =2$ then $Q(-1)$ has 
 $c_1(Q(-1))=0,\, c_2(Q(-1)) =1. $  
 Now  Riemann-Roch theorem together with Proposition (7.1) and Theorem
 (7.4) of \cite{Ha1} we see that $h^0(Q(-1)) = 1.$
 Hence 
 there is an exact sequence 
 \begin{equation}\label{eq1}
 0 \to \OOO_{\PP^2}(1) \to Q \to I_p(1) \to 0,
 \end{equation}
 where $I_p$ is the ideal sheaf of a point $p \in \PP^2.$
From the exact sequence \eqref{eq1} we see that $h^0(Q) = 5$ and $Q$ 
is generated by sections. Since $Q$ is a rank two vector bundle  on 
$\PP^2$ generated by section we see that $Q$ is generated by
$4$ independent sections. Hence we obtain a surjective 
morphism  of bundles
$$ \OOO_{\PP^2}^4 \to Q \to 0$$
which determines a morphism $\phi : \PP^2 \to Gr(2, \CC^{4}) $ as 
required.  $\hfill{\Box}$

Let $Q$ be rank two vector bundle on $\PP^2$  with 
$c_1(Q) =2,\, c_2(Q) =2.$ By Lemma (\ref{lem1}) $Q$ determines a
morphism $\phi : \PP^2 \to Gr(2, \CC^{4}). $ 

\begin{lem}\label{lem2}
Let $Q$ be a semi-stable rank two vector bundle 
on $\PP^2$ with $c_1(Q) =2,\, c_2(Q) =2.$ 
Then there is  an exact sequence
\begin{equation}\label{eq2}
0\to \OOO_{\PP^2}(-1) \to \OOO_{\PP^2}^2\oplus \OOO_{\PP^2}(1) \to Q \to 0
\end{equation}
of vector bundles on $\PP^2.$
\end{lem}
 
{\bf Proof:}  From the proof of Lemma (\ref{lem1}) it follows that
 $Q$ fits into an exact sequence \eqref{eq1} and the induced map
 $$ {\rm H}^0(\PP^2, Q) \to {\rm H}^0(\PP^2, I_p(1))$$ 
 of cohomology groups is surjective.
Observe that if\  $p  \in \PP^2$ then\\
${\rm H}^0(\PP^2, I_p(1))$ 
is two dimensional vector space and the natural map
 $$ {\rm H}^0(\PP^2, I_p(1))\otimes \OOO_{\PP^2} \to I_p(1) $$
 is surjective. Thus from \eqref{eq1} and the above observations 
 we see that there is a surjection of vector bundles
 $$ \OOO_{\PP^2}(1)\oplus \OOO_{\PP^2}^2 \to Q \to 0. $$ 
 A simple Chern class computation will show that the kernel of this
 surjection is equal to the line bundle $\OOO_{\PP^2}(-1)$ and hence 
 we get the existence of required exact 
 sequence \eqref{eq2}.  $\hfill{\Box}$
 
\section{Explicit constructions}

 If $X,Y$ and $Z$ is the standard basis of ${\rm H}^0(\PP^2, \OOO_{\PP^2}(1)),$
then
the global section $(X,Y, Z^2)$ of the vector bundle 
$\OOO_{\PP^2}(1)^2\oplus \OOO_{\PP^2}(2)$ is nowhere vanishing on
$\PP^2.$ Thus the bundle map
$$ 0\to \OOO_{\PP^2}(-1) \to \OOO_{\PP^2}^2\oplus \OOO_{\PP^2}(1)$$
given by  $(X,Y, Z^2)$ is injective. If $Q$ is the cokernel  of this 
injection then   
\begin{equation}\label{eq3}
0\to \OOO_{\PP^2}(-1){\stackrel{(X,Y, Z^2)}{\longrightarrow}} \OOO_{\PP^2}^2\oplus \OOO_{\PP^2}(1) \to Q \to 0
\end{equation}
is an exact sequence of vector bundles on $\PP^2$  with  rank of $Q$ 
two and  $c_1(Q) =2,\, c_2(Q) =2.$
From cohomology exact sequence associated to \eqref{eq3} we see that
\begin{equation}\label{eq4}
{\rm H}^0(\PP^2, \OOO_{\PP^2}^2\oplus \OOO_{\PP^2}(1) ) 
\simeq {\rm H}^0(\PP^2 , Q).
\end{equation} 

 The  vector space ${\rm H}^0(\PP^2, \OOO_{\PP^2}^2\oplus \OOO_{\PP^2}(1) )$
is equal to 
$$ \oplus_{i=1}^5 \CC v_i,$$
where $v_1=(1,0,0), v_2=(0,1,0), v_3=(0,0,X), v_4=(0,0,Y), 
v_5=(0,0,Z).$ If  $w_i \in {\rm H}^0(\PP^2 , Q)$ 
is the image of $v_i$ under the isomorphism of \eqref{eq4}
then $w_i, i=1,\ldots,5$ is a basis of ${\rm H}^0(\PP^2 , Q).$

\begin{lem}\label{lem3}
If $E$ is vector bundle of rank $2$ on $\PP^2$ generated by three global
sections then $c_1(E) = d$ and $c_2(E) =d^2$ for some integer $d\geq 0.$
\end{lem}
{\bf Proof:} If a rank two vector bundle $E$ is generated by three
global sections then we get an exact sequence
$$ 0\to \OOO_{\PP^2}(-d) \to \OOO_{\PP^2}^3 \to E \to 0 $$
of vector bundles on $\PP^2$ and hence the required result. $\hfill{\Box}$

From Lemma(\ref{lem3}) we see that the bundle $Q$ on $\PP^2$ that
we have constructed above cannot be generated by $3$ sections.

{\bf Example 1}: Let $Q$ be the vector bundle of rank two on $\PP^2$
defined by the exact sequence \eqref{eq3}. If $w_i, 1\leq i \leq 5$
are  the sections of $Q$ defined above, then the set 
$S_1= \{w_i; i=1,\ldots,4. \}$ is a generating set of sections of $Q,$ 
i.e.,if  $S=\{e_i; i=1,\ldots,4.\}$ is the standard basis of $\OOO_{\PP^2}^4$ 
then the bundle map
$$ \OOO_{\PP^2}^4 \to Q$$ 
obtained by sending $e_i$ to $w_i$ for $i=1,\ldots, 4$ is surjective. 
Thus we get a morphism 
$$ \phi_{S_1} : \PP^2 \to Gr(2, \CC^{4}). $$
If $ p: Gr(2, \CC^{4}) \to \PP^5 $ is the Plucker imbedding, then
the morphism 
$$ p\circ\phi_{S_1} : \PP^2 \to \PP^5 $$
is given by $(x,y,z) \mapsto (z^2, -xy, -y^2, x^2, xy, 0).$
Let $z_i,(0\leq i \leq 5)$ be the homogeneous coordinates of $\PP^5.$ 
Then the image of $\PP^2$ under this morphism is a rank 3 
quadric $V$ in a linear 
$\PP^3(\subset \PP^5)$: In fact $V=Z(Z_5, Z_1+Z_4, Z_1^2+Z_2Z_3)).$
 Let $p=(1,0,0,0,0,0) \in V$ and 
$C = V \cap H$ where $H$ is the hyperplane of $\PP^5$ defined by
$Z_0.$ The morphism  
$$\phi_{S_1}|_{\PP^2-\{(0,0,1)\}} =: \PP^2-\{(0,0,1)\} \to V-\{p\}$$
is a two sheeted ramified covering 
ramified precisely  along $C.$ Moreover $\phi_{S_1}((0,0,1)) = p$ and
the differential map  
${d\phi_{S_1}}_{(0,0,1)}$ is zero. In this case we see that the image of
$\PP^2$ in $Gr(2, \CC^{4})$ is a singular surface with exactly one
singularity.

{\bf Example 2}: Let $Q$ be the vector bundle of rank two on $\PP^2$
defined by the exact sequence \eqref{eq3} and let $w_i, 1\leq i \leq 5$
be the sections of $Q$ defined above and let $u_1= w_1, 
u_2= w_2, u_3= w_3, 
u_4=w_4-w_5.$ The set 
$S_2= \{u_i; i=1,\ldots,4 \}$ is a generating set of sections of $Q,$ 
i.e.,if  $S=\{e_i; i=1,\ldots,4.\}$ is the standard basis of 
$\OOO_{\PP^2}^4$ then the bundle map
$$ \OOO_{\PP^2}^4 \to Q$$ 
obtained by sending $e_i$ to $u_i$ for $i=1,\ldots, 4,$ is surjectve. 
Thus we get a morphism 
$$ \phi_{S_2} : \PP^2 \to Gr(2, \CC^{4}). $$
If $ p: Gr(2, \CC^{4}) \to \PP^5 $ is the Plucker imbedding, then
we  see that the morphism 
$$ p\circ\phi_{S_2} : \PP^2 \to \PP^5 $$
is given by 
$$(x,y,z)\mapsto (z^2, -xy, -y^2+yz, x^2, xy-xz, 0).$$    
The image of  $\PP^2$ under this morphism is an intersection of two
quadrics $V$  in a linear $\PP^4$: In fact 
$V=V(Z_5, Z_1Z_4-Z_2Z_3, (Z_1+Z_4)^2-Z_0Z_3).$

{\bf  Example 3} :
The  vector space ${\rm H}^0(\PP^2, \OOO_{\PP^2}^2\oplus \OOO_{\PP^2}(1) )$
is equal to 
$$ \oplus_{i=1}^5 \CC v_i,$$
where $v_1=(1,0,0), v_2=(0,1,0), v_3=(0,0,X), v_4=(0,0,Y), 
v_5=(0,0,Z).$ If  $w_i \in {\rm H}^0(\PP^2 , Q)$ 
is the image of $v_i$ under the isomorphism of \eqref{eq4}
then $w_i, i=1,\ldots,5$ is a basis of ${\rm H}^0(\PP^2 , Q).$

Let $Q$ be the vector bundle of rank two on $\PP^2$
defined by the exact sequence \eqref{eq3} and let $w_i, 1\leq i \leq 5$
be the sections of $Q$ defined above and let $u_1= w_1, 
u_2= w_2, u_3= w_3+d w_4, 
u_4=a w_4+w_5,$ where $a,d$ non zero complex numbers. The set 
$S_3= \{u_i; i=1,\ldots,4 \}$ is a generating set of sections of $Q,$ 
i.e.,if  $S=\{e_i; i=1,\ldots,4.\}$ is the standard basis of 
$\OOO_{\PP^2}^4$ then the bundle map
$$ \OOO_{\PP^2}^4 \to Q$$ 
obtained by sending $e_i$ to $u_i$ for $i=1,\ldots, 4,$ is surjectve. 
Thus we get a morphism 
$$ \phi_{S_3} : \PP^2 \to Gr(2, \CC^{4}). $$
If $ p: Gr(2, \CC^{4}) \to \PP^5 $ is the Plucker imbedding, then
we  see that the morphism 
$$ p\circ\phi_{S_3} : \PP^2 \to \PP^5 $$
is given by 
$$(x,y,z)\mapsto (z^2, -(x+d y)y, -(ay+z)y, (x+d y)x, (ay+z)x, 0).$$    
The image of  $\PP^2$ under this morphism can be seen to be equal 
is equal to  intersection of two
independent singular quadrics in $\PP^4.$

\section{Main Theorem}

\begin{thm}\label{image}
 Let $\phi:\PP^2 \to Gr(2, \CC^{4})$ be a morphism. Assume that 
 $c_1(Q)=2$ and $c_2(Q)= 2,$ where $Q$ is the pull back
 by $\phi$ of the universal rank two quotient bundle on $Gr(2, \CC^{4}).$
  Then the image of $\PP^2$ in $Gr(2, \CC^{4})$ is
 \begin{itemize}
  \item[a)] either a complete intersection of two independent hyperplanes,
  \item[b)] or is a complete intersection of a hyperplane and a quadric,
 \end{itemize}
Here a hyperplane(respectively, quadric)
 means a divisor in the class of the ample generator (respectively, twice the 
class of the  ample generator)
 of the the Picard group of $Gr(2, \CC^{4}).$ In case of a) the image is 
 isomorphic to cone in $\PP^3$ over a conic. In the  case of b) the image 
  is isomorphic to base locus of a pencil consists of  singular
quadrics in $\PP^4$ of rank $3$ and $4.$ More over the image surface is 
singular exactly along a line of $\PP^4$ contained in the surface.
\end{thm}

{\bf Proof:}  By our assumption, the vector bundle $Q$ on $\PP^2$
is generated by global sections. By Lemma \ref{lem2},
$Q$ fits into an exact sequence \eqref{eq2}. In the equation \eqref{eq2} the bundle map 
$$ 0\to \OOO_{\PP^2}(-1) \to \OOO_{\PP^2}^2\oplus \OOO_{\PP^2}(1) $$
is given by $s=(A, B, Q),$ where $A, B\in \text{H}^0(\OOO_{\PP^2}(1))$ 
and $Q \in \text{H}^0(\OOO_{\PP^2}(2))$
with out common zeros in $\PP^2.$   Then
the morphism $\phi$ is determined by
four linearly independent global sections of the bundle 
$E=\OOO_{\PP^2}^2\oplus \OOO_{\PP^2}(1)$ whose images generate the 
bundle $Q.$  If $w_1,w_2,w_3,w_4 $ four linearly independent global 
sections of the bundle 
$E=\OOO_{\PP^2}^2\oplus \OOO_{\PP^2}(1)$ whose images generate the 
bundle $Q$ and $u_1,u_2,u_3,u_4 $ be any other basis of the vector
space generated by $w_1,w_2,w_3,w_4 ,$ then the morphism 
defined by $w_1,w_2,w_3,w_4 $ and $u_1,u_2,u_3,u_4 $ differ by an
automorphism of $Gr(2, \CC^{4}).$ Let $w_i=(a_i,b_i, g_i), \,
1\leq i \leq 4$ be four global sections of 
$E=\OOO_{\PP^2}^2\oplus \OOO_{\PP^2}(1)$ such that their images 
in $Q$ generate $Q.$ Since $Q$ is not direct sum of line bundles
we see that sections of the bundle $\OOO_{\PP^2}^2$ given by 
$(a_i, b_i), \, 1\leq i \leq 4$ generate $\OOO_{\PP^2}^2.$ Hence  
by taking suitable linear combinations we can assume 
that the four linearly independent global sections 
are of the form
$w_1=(1,0,f_1), w_2=(0,1,f_2), w_3=(0,0,f_3), w_4=(0,0,f_4),$ where
$f_3, f_4 \in \text{H}^0(\OOO_{\PP^2}(1))$ are linearly independent.
If $ p: Gr(2, \CC^{4}) \to \PP^5 $ is the Plucker imbedding, then
we  see that the morphism 
$$ p\circ\phi : \PP^2 \to \PP^5 $$
is given by 
$$(x;y;z)\mapsto$$
$$(D_0(x,y,z);D_1(x,y,z); D_2(x,y,z); D_3(x,y,z);
D_4(x,y,z); 0),$$  
where $D_0=Q -f_1A, D_1=-(f_3B),D_2=-(f_4B),D_3=-f_3A, D_4=-f_4B.$

{\bf Case 1:} The subspace $<f_3, f_4> $ generated by $f_3, f_4$ is equal to
$<A, B>.$ In this case the four quadrics $-f_3B, -f_4B, f_3A, f_4A$ on 
$\PP^2$ are linearly dependent, say $a (-f_3B)+b( -f_4B)+c (f_3A)+d( f_4A)=0,$
with $(a,b,c,d)\neq (0,0,0,0). $
Then  the image of $\PP^2$ in $\PP^5$ is an intersection of two independent hyperplanes namely, $Z_5= 0$
and $a Z_1+b Z_2+c Z_3+d Z_4=0,$ with the Plucker embedding of the 
Grassmannnian given by $Z_0Z_5+Z_1Z_4-Z_2Z_3,$ where $Z_0, \ldots ,Z_5$ are the Plucker
coordinate functions. Hence the image of $\PP^2$ in $ Gr(2, \CC^{4})$
 is as stated in the case a) of the theorem.

{\bf Case 2:} The subspace $<f_3, f_4> $ generated by $f_3, f_4$ is not equal to
$<A, B>.$  In this case, if there is a non-trivial linear relation 
$$a(Q -f_1B)+b(-f_3B)+c (-f_4B)+d (f_3A)+e (f_4A) = 0 $$
exists, then again we see that the image of $\PP^2$ in $\PP^5$ under
$p\circ\phi $ is an intersection of two independent hyperplanes namely, 
$Z_5= 0$
and $a Z_0+b Z_1+c Z_2+d Z_3+e Z_4=0$ with the Plucker embedding of the 
Grassmannnian given by $Z_0Z_5+Z_1Z_4-Z_2Z_3.$ Hence the image 
of $\PP^2$ in $ Gr(2, \CC^{4})$
 is as stated in the case a) of the theorem.
 
 Continuing the proof in the case 2, we can assume image is not as in a) in which case
 degree of the image of $\PP^2$ in $\PP^5$ under $p\circ\phi $ 
 has to be four. As 
 above set
 $$D_1=(Q -f_1B), D_2=(-f_3B), D_3= (-f_4B), D_4=(f_3A), D_5=(f_4A) , $$
 then we see that
 \begin{equation}\label{eq5}
  D_2D_5-D_3D_4 = 0.
 \end{equation}
 
 This shows that image of $\PP^2$  under $p\circ\phi $ is a degree
 four surface contained in intersection of the hyperplane $Z_5=0$
 with the plucker embedding of the Grassmannian 
 given by the equation $Z_0Z_5+Z_1Z_4-Z_2Z_3 =0,$ where $Z_i\, (0\leq
 i \leq 5)$ are the homogeneous coordinate functions on $\PP^5.$
 By Tango's result  about the embeddings of $\PP^2$ in 
 Grassmannian (see, \cite{Ta}) we 
 conclude that the image of $\PP^2$ under  $p\circ\phi $ is 
 singular surface of degree four in $\PP^4$ and is given by 
 base point free linear system of conics. Since the linear system is of dimension  four we see that the image of $\PP^2$ under  $p\circ\phi $
 is projection of the Veronese surface $V \subset \PP^5$ from a point
 $p \in sec(V)\subset \PP^5$ and $p \notin V,$ where $sec(V)$ denotes
 the secant variety of $V.$ Hence by Remark (\ref{point}) it follows that 
 the image of $\PP^2$ in $\PP^4$ given by $p\circ\phi $ is cut out by a linear pencil of quadrics in $\PP^4.$
 Thus the image of $\PP^2$ under $\phi$ in this case is a complete intersection of a hyperplane and a quadric as stated in b).
 
The last statement about the singularities can be easily checked.
$\hfill{\Box}$

{\it Acknowledgements}: We thank referee for his valuable suggestions.
Last named author would like to thank University
of Lille-1 at Lille and University of Artois at Lens. He would also like to 
thank university of Paris 6 and IRSES-Moduli program.

\end{document}